\documentclass[letterpaper, 10 pt, conference]{ieeeconf}  

\IEEEoverridecommandlockouts                              
                                                          
\usepackage{xspace,amssymb,epsfig,subfigure,syntonly,cite}
\usepackage{epsfig,amsmath,color}
\usepackage{xspace,syntonly}
\usepackage{url}
\usepackage{algorithm,algorithmic}

\newcommand{\utwi}[1]{\mbox{\boldmath $#1$}}

\newcommand{\trace}{{\textrm{Tr}}}
\newcommand{\rank}{{\textrm{rank}}}

\newcommand{\cD}{{\cal D}}

\newcommand{\cN}{{\cal N}}

\newcommand{\cE}{{\cal E}}

\newcommand{\cU}{{\cal U}}

\newcommand{\cV}{{\cal V}}

\newcommand{\cY}{{\cal Y}}

\newcommand{\bb}{{\bf b}}
\newcommand{\bd}{{\bf d}}
\newcommand{\be}{{\bf e}}
\newcommand{\bbf}{{\bf f}}
\newcommand{\bg}{{\bf g}}
\newcommand{\bh}{{\bf h}}

\newcommand{\br}{{\bf r}}
\newcommand{\bs}{{\bf s}}
\newcommand{\bx}{{\bf x}}
\newcommand{\bu}{{\bf u}}
\newcommand{\bw}{{\bf w}}
\newcommand{\bv}{{\bf v}}
\newcommand{\bi}{{\bf i}}
\newcommand{\bz}{{\bf z}}
\newcommand{\by}{{\bf y}}
\newcommand{\bA}{{\bf A}}

\newcommand{\bC}{{\bf C}}
\newcommand{\bD}{{\bf D}}

\newcommand{\bI}{{\bf I}}
\newcommand{\bX}{{\bf X}}

\newcommand{\bY}{{\bf Y}}
\newcommand{\bV}{{\bf V}}

\newcommand{\blambda}{{\utwi{\lambda}}}

\newcommand{\bPsi}{{\utwi{\Psi}}}

\newcommand{\bPhi}{{\utwi{\Phi}}}

\newcommand{\bxi}{{\utwi{\xi}}}

\newcommand{\bUpsilon}{{\utwi{\Upsilon}}}

\newcommand{\sfH}{\textsf{H}}
\newcommand{\sfT}{\textsf{T}}

\usepackage{epstopdf}

\usepackage{color}

\begin{document}

\newtheorem{definition}{Definition}
\newtheorem{remark}{Remark}
\newtheorem{proposition}{Proposition}
\newtheorem{lemma}{Lemma}
\newtheorem{assumption}{Assumption}
\newtheorem{theorem}{Theorem}



\title{\LARGE \bf
Regulation of Dynamical Systems to Optimal Solutions of Semidefinite Programs: Algorithms and Applications to  AC Optimal Power Flow
}

\author{Emiliano Dall'Anese, Sairaj V. Dhople, and Georgios B. Giannakis
\thanks{The Authors were supported in part by the National Science Foundation (NSF) through grants NSF-CCF grant no. 1423316 and CyberSEES grant no. 1442686. S. V. Dhople was also supported in part by the NSF CAREER award ECCS-1453921.}
\thanks{E. Dall'Anese is with the National Renewable Energy Laboratory, Golden, CO, USA. S. V. Dhople and G. B. Giannakis are with the Dept. of ECE and Digital Tech. Center, University of Minnesota, 
Minneapolis, MN, USA.   E-mail: {emiliano.dallanese@nrel.gov}, {\{sdhople, georgios\}@umn.edu}%
}
}

\maketitle

\begin{abstract}

This paper considers a collection of networked nonlinear dynamical systems, and addresses the synthesis of feedback controllers that seek optimal operating points corresponding to the solution of pertinent network-wide optimization problems. Particular emphasis is placed on the solution of semidefinite programs (SDPs). The design of the feedback controller is grounded on a dual $\epsilon$-subgradient approach, with the dual iterates utilized to dynamically update the dynamical-system reference signals. Global convergence is guaranteed for diminishing stepsize rules, even when the reference inputs are updated at a faster rate than the dynamical-system settling time. The application of the proposed framework to the control of power-electronic inverters in AC distribution systems is discussed. The objective is to bridge the time-scale separation between real-time inverter control and network-wide optimization. Optimization objectives assume the form of SDP relaxations of prototypical AC optimal power flow problems.
\end{abstract}

\section{Introduction}
\label{sec:Introduction}

This paper addresses the synthesis of feedback controllers that seek to regulate networked nonlinear dynamical systems to the optimal solution of a convex constrained optimization problem. The setup is relevant in several multi-agent-system applications. In this context, time-scale separation is typically leveraged to enforce a strict temporal barrier in terms of when the optimal setpoints are solved for and dispatched to the dynamical systems. However, operational efficiency can be improved by compressing the time scales and devising means to synergize the implementation of the optimization problems and real-time controllers. 

Previous efforts in this domain are grounded in the seminal work~\cite{Arrowbook58}, where dynamical systems serve as proxies for optimization variables and multipliers, and are synthesized to evolve in a gradient-like fashion to the saddle points of the Lagrangian function associated with the convex optimization problem~\cite{DeHaan04,Elia-CDC11,Brunner-CDC12}. 
Particularly relevant to this paper are the results reported in~\cite{Jokic09}, where  a continuous-time feedback controller that seeks Karush-Kuhn-Tucker (KKT) conditions for optimality of a convex constrained optimization problem is developed. A heuristic comprising continuous-time dual ascent and discrete-time reference-signal updates is considered in~\cite{Hirata}, and local stability of the resultant closed-loop system is established. Distinct from~\cite{Jokic09,Hirata} as well as previous efforts in e.g.,~\cite{DeHaan04,Elia-CDC11,Brunner-CDC12,Elia-Allerton13,NaLi_ACC14,Chen_CDC14,Jokic09,Hirata}, this work leverages dual-subgradient methods to develop a feedback controller that steers the dynamical-system outputs towards the solution of a  convex constrained optimization problem. The proposed scheme involves the update of dual and primal variables in a discrete-time fashion, with the latter constituting the reference-input signals for the dynamical systems. When dual and primal variables are updated at a faster rate than the system settling time, it is shown that the dual ascent step is in fact an $\epsilon$-subgradient~\cite{Kiwiel04}. This is particularly relevant in settings where the reference signals may be updated continuously (within the limits of affordable computational burden), without necessarily waiting for the underlying dynamical systems to converge to intermediate reference levels. Convergence of system outputs to the solution of SDP-type problems is established with diminishing stepsize rules and strictly convex cost functions. Although the framework is outlined for a semidefinite program (SDP), similar convergence claims can be established for other types of optimization problems.  

The application of the proposed framework in the context of power systems is discussed, with particular emphasis on distribution networks featuring power-electronic-inverter-interfaced (renewable) energy resources~\cite{Iravanibook10}. In particular, the controller devised in this paper is utilized to steer the output of inverters towards the solution of an AC optimal power flow (OPF) problem, which yields steady-state active- and reactive-power injections that are optimal according to well-defined optimization criteria. Since the AC OPF task corresponds to a nonconvex optimization problem, an SDP relaxation~\cite{LavaeiLow,Dallanese-TSG13} is leveraged. In this context, the objective is to bridge the temporal gap between long-term energy management and real-time control, to ensure adaptability to changing ambient conditions and loads, and guarantee seamless renewable energy integration without compromising system stability~\cite{OID,OID_TEC,Dorfler14}. Similar controllers focused on an economic dispatch problem have been proposed for bulk power systems in~\cite{Jokic_JEPES}.  Modified automatic generation and frequency control methods that incorporate optimization objectives corresponding to DC optimal power flow (OPF) problems are proposed for lossless bulk power systems in~\cite{NaLi_ACC14,Chen_CDC14}.  Strategies that integrate economic optimization within droop control for islanded lossless microgrids are developed in~\cite{Dorfler14}. Different from the continuous-time controllers developed in~\cite{NaLi_ACC14,Chen_CDC14, Jokic_JEPES},  the proposed approach accounts for computational limits in the update of the inverter setpoints (which naturally lead to discrete-time reference updates), and considers strict inverter-generation limits.

\section{Problem Formulation}
\label{sec:ProblemFormulation}
Consider $N$ nonlinear dynamical systems described by\footnote{\emph{Notation}. Upper-case (lower-case) boldface
letters will be used for matrices (column vectors); $(\cdot)^\sfT$ for transposition; $(\cdot)^*$ complex-conjugate; and, $(\cdot)^\sfH$ complex-conjugate transposition;  $\Re\{\cdot\}$ and $\Im\{\cdot\}$ denote the real and imaginary parts of a complex number, respectively; $\mathrm{j} := \sqrt{-1}$. $\trace(\cdot)$ the matrix trace; $\mathrm{rank}(\cdot)$ the matrix rank; $|\cdot|$ denotes the magnitude of a number or the cardinality of a set; $\mathrm{vec}(\bX)$ returns a vector stacking the columns of matrix $\bX$, and $\mathrm{bdiag}(\{\bX_i\})$ forms a block-diagonal matrix. $\mathbb{R}^N$ and $\mathbb{C}^N$ denote the spaces of $N\times1$ real-valued and complex-valued vectors, respectively; $\mathbb{N}$ the set of natural numbers; and, $\mathbb{H}_+^{N \times N}$ denotes the space of $N \times N$ positive semidefinite Hermitian matrices. Given vector $\bx$ and square matrix $\bX$,  $\|\bx\|_2$ denotes the Euclidean norm of $\bx$, and $\|\bV\|_2$ the (induced) spectral norm of matrix $\bX$. $[\bx]_i$ ($[\bf(\bx)]_i$) points to the $i$-th element of a vector $\bx$ (vector-valued function $\bbf(\bx)$). $\dot{\bx}(t)$ is the time derivative of $\bx(t)$. Given a scalar function $f(\bx): \mathbb{R}^n \rightarrow \mathbb{R}$, $\nabla_\bx f(\bx)$ returns the gradient $[\frac{\partial f}{\partial x_1}, \ldots, \frac{\partial f}{\partial x_n}]^\sfT$. For a continuous function $f(t)$, $f[t_k]$ denotes its value sampled at $t_k$. Finally, $\bI_N$ denotes the $N \times N$ identity matrix; and, $\mathbf{0}_{M\times N}$, $\mathbf{1}_{M\times N}$ the $M \times N$ matrices with all zeroes and ones, respectively.}  
\begin{subequations}
\label{eq:sys1}
\begin{align}
\dot{\bx}_i(t) & = \bbf_i \Big( \bx_i(t), \bd_i(t), \bu_i(t) \Big) \label{eq:sys1-st} \\
\by_i(t) & =  \br_i \Big( \bx_i(t), \bd_i(t) \Big) , \quad i \in \mathcal{N} := \{1, \ldots, N \} \label{eq:sys1-obs} 
\end{align}
\end{subequations}
where: $\bx_i(t) \in \mathbb{R}^{n_{x,i}}$ is the state of the $i$-th dynamical system at time $t$; $\by_i(t) \in \cY_i \subset \mathbb{R}^{n_{y,i}}$ is the measurement of state $\bx_i(t)$ at time $t$; $\bu_i(t) \in \cY_i$ is the reference input; and $\bd_i(t) \in \cD_i  \subset \mathbb{R}^{n_{d,i}}$ is the exogenous input. Finally, $\bbf_i: \mathbb{R}^{n_{x,i}} \times \mathbb{R}^{n_{z,i}} \times \mathbb{R}^{n_{d,i}} \times \mathbb{R}^{n_{y,i}} \rightarrow \mathbb{R}^{n_{x,i}}$ and $\br_i:  \mathbb{R}^{n_{x,i}} \times \mathbb{R}^{n_{d,i}}  \rightarrow  \mathbb{R}^{n_{y,i}}$ are arbitrary (non)linear functions. Similar to, e.g.,~\cite{Jokic09,Hirata}, the following system behavior for constant exogenous inputs and reference signals is presumed.

\vspace{.1cm}

\begin{assumption}
\label{ass:DynSystems}
For given \emph{constant} exogenous inputs $\{\bd_i \in \cD_i\}_{i \in \cN}$ and reference signals $\{\bu_i \in \cY_i\}_{i = 1}^N$, there exist equilibrium points $\{\bx_i\}_{i = 1}^N$ for~\eqref{eq:sys1} that satisfy: 
\begin{subequations}
 \label{eq:sys1-equilibrium}
\begin{align}
\mathbf{0} & = \bbf_i \left(\bx_i, \bd_i, \bu_i \right) \label{eq:sys1-st-equilibrium} \\
\bu_i & =  \br_i \left( \bx_i, \bd_i \right), \quad i \in \cN  \, . \label{eq:sys1-obs-equilibrium} 
\end{align}
Notice that in~\eqref{eq:sys1-obs-equilibrium} the equilibrium output coincides with the commanded input $\bu_i$. 
Furthermore, these equilibrium points are assumed to be locally asymptotically stable.  \hfill $\Box$ 
\end{subequations}
\end{assumption}

\vspace{.1cm}

For given exogenous inputs $\{\bd_i \in \cD_i\}_{i = 1}^N$, consider the following optimization problem associated with $\{\bu_i\}_{i = 1}^N$: 
\begin{subequations}
\label{eq:p1}
\begin{align}
\mathrm{(P1)} \hspace{0cm}  \min_{\bV \in \cV, \{\bu_i \in \cY_i \}}& \,\, H(\bV) + \sum_{i \in \cN} \left(  \frac{1}{2}\bu_i^\sfT \bA_i \bu_i + \bb_i^\sfT \bu_i \right) \label{eq:p1-cost} \\
& \hspace{-1.8cm} \mathrm{subject~to}  \,\, \bh_i(\bV) + \bg_i(\bu_i,\bd_i)  = \mathbf{0} , \, \forall \, i \in \cN \label{eq:p1-eq} 
\end{align}
\end{subequations}
where $\cV \subset \mathbb{H}_+^{n_V \times n_V}$  is a convex, closed, and bounded subset of the cone of positive semidefinite (Hermitian) matrices; function $H(\bV): \mathbb{H}_+^{n_V \times n_V} \rightarrow \mathbb{R}$ is known, strictly convex and finite over $\cV$; $\bA_i \succ \mathbf{0}$ and $\bb_i \in \mathbb{R}^{n_{x,i}}, \forall i \in \cN_D$; the vector-valued function $\bh_i(\bV): \mathbb{H}_+^{n_V \times n_V} \rightarrow \mathbb{R}^{n_{y,i}}$ is affine; and, $\bg_i(\bu_i, \bd_i): \mathbb{R}^{n_{y,i}} \times \mathbb{R}^{n_{d,i}} \rightarrow \mathbb{R}^{n_{y,i}}$  takes the form $\bg_i(\bu_i, \bd_i) = \bC_i \bu_i + \bD_i \bd_i$, 
with $\bC_i \in \mathbb{R}^{n_{y,i} \times n_{y,i}}$ and $\bD_i \in \mathbb{R}^{n_{y,i} \times n_{d,i}}$ known. Finally, sets $\{\cY_i\}_{i \in \cN_D}$, which define the space of possible reference inputs for the dynamical systems, are assumed to comply to the following requirement.

 \vspace{.1cm}
  
\begin{assumption}
\label{ass:compactSetY}
Sets $\{\cY_i\}_{i = 1}^N$ are convex, closed, and bounded. \hfill $\Box$ 
\end{assumption}

\vspace{.1cm}
 
\noindent With these assumptions, problem $\mathrm{(P1)}$ is a \emph{convex} program; moreover, it can be reformulated into a standard SDP form by resorting to the epigraph form of the cost function. 

It is evident from~\eqref{eq:sys1-obs-equilibrium} that $\mathrm{(P1)}$ defines the optimal operating setpoints of the dynamical systems~\eqref{eq:sys1} in terms of steady-state outputs~\cite{Jokic09,Hirata}. In fact, by utilizing the optimal solution $\{\bu_i^{\textrm{opt}}\}_{i \in \cN_D}$ of $\mathrm{(P1)}$ as reference inputs, it follows from~\eqref{eq:sys1-obs-equilibrium} that each system output will eventually be driven to the point $\by_i = \bu_i^{\textrm{opt}}$. 

In principle, $\mathrm{(P1)}$ could be solved centrally by a system-level control unit~\cite{LavaeiLow} or in a decentralized fashion~\cite{OID_TEC,Dorfler14}, and the  reference signals $\{\bu_i^{\textrm{opt}}\}_{i \in \cN_D}$ could be subsequently dispatched for the dynamical systems. In lieu of this solution with strict temporal boundaries, the \emph{objective} here is to design a decentralized feedback controller for the dynamical systems~\eqref{eq:sys1}, so that the resultant closed-loop system is globally convergent to an equilibrium point $\{\bx_i\}_{i = 1}^N$, $\{\by_i =  \br_i ( \bx_i, \bd_i )\}_{i = 1}^N$, where the values $\{\by_i\}_{i = 1}^N$ of the steady-state outputs coincide with the optimal solution $\{\bu_i^{\textrm{opt}}\}_{i = 1}^N$ of $\mathrm{(P1)}$.

\section{Feedback Controller Synthesis}
\label{sec:Subgradient}

\subsection{A Primer on Dual Gradient Methods}
\label{sec:Subgradient_prel}

To streamline exposition, it is  convenient to consider expressing the linear equality constraints~\eqref{eq:p1-eq} in the compact form $\bh(\bV) + \bg(\bu,\bd) = \mathbf{0}$, where $\bu := [\bu_1^\sfT, \ldots, \bu_N^\sfT]^\sfT$,  $\bd := [\bd_1^\sfT, \ldots, \bd_N^\sfT]^\sfT$, $\bh(\bV) := [\bh_1^\sfT(\bV), \ldots, \bh_N^\sfT(\bV)]^\sfT$, and 
\begin{align}
\bg(\bu, \bd) := \bC \bu + \bD \bd , \label{eq:gunction_g2}
\end{align}
with $\bC$ denoting the block-diagonal matrix specified as $\bC := \mathrm{bdiag}(\{\bC_i\}_{i = 1}^N)$ and $\bD$ formed using $\{\bD_i\}_{i = 1}^N$. Thus, recalling that $[\bx]_i$ ($[f(\bx)]_i$) denotes the $i$-th element of a vector $\bx$ (vector-valued function $f(\bx)$), the  following is assumed for the convex program $\mathrm{(P1)}$.

\vspace{.1cm}

\begin{assumption}
\label{ass:constraintQualification}
Problem $\mathrm{(P1)}$ has a non-empty feasible set and a finite optimal cost. Furthermore, the vectors 
\begin{align}
\nabla_{[\textrm{vec}^\sfT (\bV), \bu^\sfT ]^\sfT} [\bh(\bV) + \bg(\bu,\bd)]_j ,\quad j = 1, \ldots, \sum_{i} n_{y,i} \label{eq:cq}
\end{align}
are linearly independent. \hfill $\Box$ 
\end{assumption}

\vspace{.1cm}

From the non-emptiness and compactness of the feasible set, and the continuity of the objective function, it follows that an optimal solution to $\mathrm{(P1)}$ exists~\cite{Elia-CDC11,Bertsekas_ConvexAnalysis}. In par with  the linear independence constraint qualification, \emph{Assumption~\ref{ass:constraintQualification}} ensures existence and uniqueness of the optimal multipliers~\cite{Wachsmuth13}. When a set of inequality constraints is added to $\mathrm{(P1)}$,~\emph{Assumption~\ref{ass:constraintQualification}} can be replaced by the Mangasarian-Fromovitz constraint qualification to ensure non-emptiness  and boundedness of the optimal multiplier set~\cite{Wachsmuth13}.

Let $\blambda_i \in \mathbb{R}^{n_{y,i}}$ denote the Lagrange multiplier associated with equality~\eqref{eq:p1-eq}, and consider the Lagrangian function corresponding to~\eqref{eq:p1}, which is defined as:  
\begin{align}
L\left(\bV, \{\bu_i\}, \{\blambda_i\}\right) &:=  \,H(\bV) + \sum_{i = 1}^N \left( \frac{1}{2} \bu_i^\sfT \bA_i \bu_i + \bb_i^\sfT \bu_i \right) \nonumber \\
& + \sum_{i = 1}^N \blambda_i^\sfT \left( \bh_i(\bV) + \bg_i(\bu_i, \bd_i) \right) \, . \label{eq:lagrangian}
\end{align}
Based on~\eqref{eq:lagrangian}, the dual function and the dual problem are given by (see, e.g.,~\cite{Bertsekas_ConvexAnalysis})
\begin{align}
q(\blambda) := \min_{\bV \in \cV, \{\bu_i \in \cY_i\}} L(\bV, \{\bu_i\},\blambda)  \label{eq:dualFunction}
\end{align}
\vspace{-.2cm}
\begin{align}
q^{\mathrm{opt}} := \max_{\blambda} \,\, q(\blambda)   \label{eq:optimalDualFunction}
\end{align}
where $\blambda := [\blambda_1^\sfT, \ldots, \blambda_N^\sfT]^\sfT$. Under current modeling assumptions, it follows that the duality gap is zero~\cite{Bertsekas_ConvexAnalysis}; furthermore, the dual function $q(\blambda)$ is concave and differentiable~\cite{Cheng87}. 

Consider utilizing a gradient method to solve the dual problem, which amounts to iteratively performing~\cite{Cheng87}:  
\begin{subequations} 
\label{eq:dualsubgradient_basic}
\begin{align}
\hspace{-.2cm} \{\bV[k], \{\bu_i[k]\}_{i = 1}^N\} & \nonumber \\
& \hspace{-1.45cm} = \arg \min_{\bV \in \cV, \{\bu_i \in \cY_i\}} \, L(\bV, \{\bu_i\}, \{\blambda_i[k]\}) \label{eq:primal_basic} \\
&\hspace{-2.8cm}  \blambda_i[k+1]  = \blambda_i[k] + \alpha_{k+1} \nabla_{\blambda_i}  L(\bV[k], \{\bu_i[k]\}, \{\blambda_i\})\,  \label{eq:dual_ascent_basic}
\end{align}
\end{subequations}
where $k \in \mathbb{N}$ denotes the iteration index,  $\alpha_{k+1} \geq 0$ is the stepsize, and~\eqref{eq:dual_ascent_basic} is repeated for all $i \in \cN$. In particular, a non-summable but square-summable stepsize sequence is adopted in this paper~\cite{Kiwiel04}; that is, there exist sequences $\{\gamma_k\}_{k \geq 0}$ and $\{\eta_k\}_{k \geq 0}$ such that:

\noindent $\mathrm{(s1)}$ $\gamma_k \rightarrow 0$ as $k \rightarrow + \infty$, and $\sum_{k = 0}^{+ \infty} \gamma_k = + \infty$;  

\noindent $\mathrm{(s2)}$ $\gamma_k \leq \alpha_k \leq \eta_k$ for all $k \geq 0$; and, 

\noindent $\mathrm{(s3)}$ $\eta_k \downarrow 0$ as $k \rightarrow + \infty$, and $\sum_{k = 0}^{+ \infty} \eta_k^2 < + \infty$.  

\noindent At iteration $k$, the same step-size $\alpha_k$ is utilized for all $i \in \cN$.  Exploiting the decomposablility of the Lagrangian, steps~\eqref{eq:dualsubgradient_basic}  can be equivalently expressed as:  
\begin{subequations}
\label{eq:dualsubgradient}
\begin{align}
& \hspace{-.2cm} \bV[k] = \arg \min_{\bV \in \cV} \, H(\bV) +  \sum_{i = 1}^N \blambda_i^\sfT[k] \, \bh_i(\bV) \label{eq:primal_V} \\
& \hspace{-.2cm} \bu_i [k] = \displaystyle{\mathrm{proj}_{\cY_i}\{- \bA_i^{-1} \bC_i^\sfT \blambda_i[k] - \bA_i^{-1} \bb_i\}} \label{eq:primal_y_2}  \hspace{-.2cm}  \\
& \hspace{-.2cm} \blambda_i[k+1] = \blambda_i[k] + \alpha_{k+1} \left( \bh_i(\bV[k]) + \bg_i(\bu_i[k], \bd_i) \right)  
 \label{eq:dual_ascent}
\end{align}
\end{subequations} 
with~\eqref{eq:primal_y_2}--\eqref{eq:dual_ascent} performed for all $i \in \cN$, and $\mathrm{proj}_{\cY}\{\bw\} := \arg \min_{\bu \in \cY} \|\bw - \bu \|_2$ denoting the projection of a vector $\bw$ onto the convex compact set $\cY$. Finally, notice that from the compactness of sets $\cV$ and $\{\cY_i\}_{i = 1}^N$, it follows that there exists a scalar $G$, $0 \leq G < + \infty$, such that
\begin{align}
\| \bh(\bV[k]) + \bg(\bu[k], \bd) \|_2 \leq G \, , \quad \forall \,\,\, k \in \mathbb{N} \, . \label{eq:bounded_gradient}
\end{align} 

Using~\eqref{eq:bounded_gradient}, and a stepsize sequence $\{\alpha_k\}_{k\geq 0}$ satisfying $\mathrm{(s1)}$--$\mathrm{(s3)}$, it turns out that the dual iterates $\blambda[k]$ converge to the optimal solution $\blambda^{\mathrm{opt}}$ of the dual problem~\eqref{eq:optimalDualFunction}; that is, $\|\blambda^{\mathrm{opt}} - \blambda[k]\|_2 \rightarrow 0$ as $k \rightarrow \infty$~\cite[Prop.~8.2.6]{Bertsekas_ConvexAnalysis},~\cite{Cheng87,Kiwiel04}. Given the strict convexity of the Lagrangian with respect to all primal variables, iterates $\bV[k]$ and $\{\bu_i[k]\}_{i = 1}^N$ become asymptotically feasible and their optimal values, $\bV^{\mathrm{opt}}$ and $\{\bu_i^{\mathrm{opt}}\}_{i = 1}^N$, can be recovered from~\eqref{eq:primal_V} and~\eqref{eq:primal_y_2}, respectively, once $\blambda^{\mathrm{opt}}$ becomes available.

\subsection{Dynamical system in-the-loop}
\label{sec:Subgradient_system}
Consider a setup where the primal and dual updates in~\eqref{eq:dualsubgradient} are performed at discrete time instants $t \in \{t_k, k \in \mathbb{N}\}$, and let $\bV[t_k]$, $\{\bu_i[t_k]\}_{i = 1}^N$, and $\{ \blambda_i[t_{k}]\}_{i = 1}^N$ denote the values of the primal and dual variables, respectively, at time $t_k$. 
With these definitions, steps~\eqref{eq:dualsubgradient} are modified to accommodate the system dynamics in~\eqref{eq:sys1} as explained next. 
At time $t_{k}$,  the system outputs are sampled as:
\begin{subequations}
\label{eq:dualsubgradient_system}
\begin{align}
\by_i[t_{k}] & =  \br_i \left( \bx_i(t_{k}), \bd_i \right) \, \,\, \forall \, i \in \cN   \label{eq:sys1-obs_sys}  
 \end{align}
and they are utilized to update the dual variables as specified in the following [cf.~\eqref{eq:dual_ascent}]:    
\begin{align}
 \blambda_i[t_{k+1}]  & = \blambda_i[t_{k}] \nonumber \\
 & \hspace{-.6cm} + \alpha_{k+1} \left( \bh_i(\bV[t_{k}]) +  \bC_i \by_i[t_{k}] + \bD_i \bd_i  \right)  \,   , \, \forall i  \, \label{eq:dual_ascent_sys}
\end{align}
where the sampled output $\by_i[t_{k}]$ replaces the primal iterate $\bu_i[t_{k}]$ on the right-hand-side of~\eqref{eq:dual_ascent_sys}. Given $\blambda_i[t_{k+1}]$,  variables $\bV[t_{k+1}]$ and $\{ \bu_i[t_{k+1}]\}_{i = 1}^{N}$ are then updated as: 
\begin{align}
& \hspace{-.1cm} \bV[t_{k+1}]  = \arg \min_{\bV \in \cV} \, H(\bV) +  \sum_{i = 1}^N \blambda_i^\sfT[t_{k+1}] \, \bh_i(\bV) \label{eq:primal_V_sys} \\
& \hspace{-.1cm} \bu_i[t_{k+1}]  = \mathrm{proj}_{\cY_i}\{ - \bA_i^{-1} \bC_i^\sfT \blambda_i[t_{k+1}] - \bA_i^{-1} \bb_i\} , \forall i .\label{eq:primal_y_sys} 
\end{align}
\end{subequations}
Once~\eqref{eq:primal_y_sys} is solved, a vector-valued reference signal taking the constant value $\bu_i[t_{k+1}]$ over $(t_{k}, t_{k+1}]$ is applied to the dynamical system~\eqref{eq:sys1-st}; i.e., $\bu_i(t) = \bu_i[t_{k+1}], t \in (t_{k}, t_{k+1}]$. At time $t_{k+1}$ the outputs $\{\by_i[t_{k+1}] \}_{i = 1}^N$ are sampled again, and~\eqref{eq:dual_ascent_sys}--\eqref{eq:primal_y_sys} are repeated. 

Steps~\eqref{eq:dual_ascent_sys}--\eqref{eq:primal_y_sys} in effect constitute the controller for the dynamical systems~\eqref{eq:sys1}. Specifically, the (continuous-time) reference signals $\{\bu_i(t)\}_{i \in \cN_D}$ produced by the controller have step changes at instants $\{t_k, k \in \mathbb{N}\}$, are left-continuous functions, and take the constant values $\{\bu_i[t_{k+1}]\}_{i \in \cN_D}$ over the time interval $(t_k, t_{k+1}]$. It is evident that if $\bu_i[t_k]$ converges to $\bu_i^{\textrm{opt}}$ as $k \rightarrow \infty$ (and thus $\bu_i(t) \rightarrow \bu_i^{\textrm{opt}}$ as $t \rightarrow \infty$), then $\by_i(t) \rightarrow \bu_i^{\textrm{opt}}$ as $t \rightarrow \infty$ by virtue of~\eqref{eq:sys1-equilibrium}. 


Suppose for now that the interval $(t_{k-1}, t_{k}]$ is large enough to allow the outputs of the dynamical systems to converge to the point $\{\bu_i[t_k]\}_{i = 1}^N$; that is, $\lim_{t \rightarrow t_{k}^-} \|\by_i(t) - \bu_i[t_k]\| = 0$, for all $k$ [cf.~\eqref{eq:sys1-equilibrium}]. In this \emph{ideal} case with a \emph{time-scale separation} between controller and system dynamics, the system dynamics do not influence the computation of the primal and dual updates, and therefore steps~\eqref{eq:dualsubgradient} and~\eqref{eq:dualsubgradient_system} coincide. The convergence results reported in Section~\ref{sec:Subgradient_prel} naturally carry over to this ideal setup.  
However, a pertinent question here is whether the closed-loop system~\eqref{eq:dualsubgradient_system} is convergent, and to what points the primal and dual iterates may converge, when at each instant $t_k$ one has that $\lim_{t \rightarrow t_{k}^-} \|\by_i(t) - \bu_i[t_k]\| \neq 0$ for at least one dynamical system; that is, no error-free tracking of the reference signals is guaranteed over each slot $(t_{k-1}, t_{k}]$.  This may represent the case where,  in an effort to compress the time scales, steps~\eqref{eq:dual_ascent_sys}--\eqref{eq:primal_y_sys} are performed continuously (within the limits of affordable computational burden), without necessarily waiting for the underlying dynamical systems to converge to the intermediate reference levels $\{\bu_i[t_k]\}_{i = 1}^N$. Or, this may represent the case where outputs are sampled without knowledge of the systems' settling times. In the following, convergence of the  closed-loop system~\eqref{eq:dualsubgradient_system} is established in this more general setup. 

For brevity, collect the system outputs in the vector $\by := [\by_1^\sfT, \ldots, \by_N^\sfT]^\sfT$. Key is to notice that, given the strict convexity of $L(\bV, \bu, \blambda[t_k])$ with respect to $\bu$, the pair $(\bV[t_k], \by[t_k])$ represents a \emph{sub-optimal} solution for the primal update~\eqref{eq:primal_basic} (and thus for~\eqref{eq:primal_V_sys}-\eqref{eq:primal_y_sys}) whenever $\lim_{t \rightarrow t_{k}^-} \|\by(t) - \bu[t_k]\| \neq 0$; that is, there exists an $\epsilon[t_k]$ such that $L(\bV[t_k], \bu[t_k], \blambda[t_k]) \leq L(\bV[t_k], \by[t_k], \blambda[t_k])$ and $L(\bV[t_k], \by[t_k], \blambda[t_k]) \leq L(\bV[t_k], \bu[t_k], \blambda[t_k]) + \epsilon[t_k]$. Thus, replacing the optimal primal iterate $\bu[t_k]$ with $\by[t_k]$ in~\eqref{eq:dual_ascent_sys} yields an $\epsilon$-subgradient step.  

Before elaborating further on the error $\epsilon[t_k]$, notice that since sets $\cV$ and $\{\cY_i\}_{i = 1}^N$ are compact, it follows that $\| \bh(\bV) + \bg(\by, \bd) \|_2$ can be bounded as [cf.~\eqref{eq:bounded_gradient}]
\begin{align}
\| \bh(\bV) + \bg(\by, \bd) \|_2 \leq G \, , \quad \forall \,\, \bV \in \cV, \, \forall \,\, \by \in \cY  \,  \label{eq:bounded_gradient_2}
\end{align} 
with $\cY := \cY_1 \times \cY_2 \times \ldots, \times \cY_N$. Furthermore, given the Lipschitz-continuity of the contraction mapping~\eqref{eq:primal_y_sys}, there exists a dual variable $\tilde{\blambda}[t_{k}]$ satisfying [cf.~\eqref{eq:primal_y_sys}]
\begin{align}
\by_i[t_{k}] = \mathrm{proj}_{\cY_i}\{- \bA_i^{-1} \bC_i^\sfT \tilde{\blambda}[t_{k}] - \bA_i^{-1} \bb_i\}, \forall i \in \cN
\label{eq:lambdatilde}
\end{align} 
that is, $\by_i[t_{k}]$ would be obtained by minimizing the Lagrangian $L(\bV, \bu, \tilde{\blambda}[t_{k}])$ when $\tilde{\blambda}[t_{k}] := [\tilde{\blambda}_1^\sfT[t_{k}], \ldots, \tilde{\blambda}_N^\sfT[t_{k}]]^\sfT$ replaces $\blambda[t_{k}]$. %
The following will be assumed for $\tilde{\blambda}[t_{k}]$.
\vspace{.1cm}

\begin{assumption}
\label{ass:bounds_error} 
There exists a scalar $\tilde{G}$, $0 \leq \tilde{G} < + \infty$, such that the bound 
\begin{align}
\| \blambda[t_k] - \tilde{\blambda}[t_k]\|_2 \leq \tilde{G}  \|\blambda[t_k] - \blambda[t_{k-1}]\|_2  \label{eq:bounded_gradient_error}
\end{align} 
holds for all $t_k$, $k \geq 1$.\footnote{Condition~\eqref{eq:bounded_gradient_error} can be re-stated in terms of the output signals $\by[t_k]$. Specifically, letting $\bxi_i[t_k] := - \bA_i^{-1} \bC_i^\sfT \blambda[t_{k}] - \bA_i^{-1} \bb_i$ be the unprojected reference signal, and assuming that matrix $\bA_i^{-1} \bC_i^\sfT$ is invertible, one has that~\eqref{eq:bounded_gradient_error} is implied by the bound $\|\bxi[t_k] - \by[t_k]\| \leq \bar{G} \|\bxi[t_k] - \bxi[t_{k-1}]\|$, upon setting $\tilde{G} = \bar{G} \|- \bA_i^{-1} \bC_i^\sfT\|_2 \|(- \bA_i^{-1} \bC_i^\sfT)^{-1}\|_2$.}  \hfill $\Box$
\end{assumption}
\vspace{.1cm}

In subsequent developments, the following bound (which originates from \emph{Assumption 4}) is leveraged: 
\begin{align}
\| \blambda[t_k] - \tilde{\blambda}[t_k]\|_2 & \leq \tilde{G}  \|\blambda[t_k] - \blambda[t_{k-1}]\|_2  \label{eq:step1} \\
&\hspace{-0.3in}  = \tilde G \| \alpha_{k} ( \bh(\bV[t_{k-1}]) + \bg(\by[t_{k-1}], \bd)) \|_2 \label{eq:step2} \\
&\hspace{-0.3in} \leq \tilde G  G \alpha_k. \label{eq:step3}
\end{align}
Note that~\eqref{eq:step2} follows from the dual update equation in~\eqref{eq:dual_ascent_sys}, and~\eqref{eq:step3} follows from~\eqref{eq:bounded_gradient_2}.

Three pertinent results that establish convergence of the overall system~\eqref{eq:dualsubgradient_system} are presented next. Lemma~\ref{eq:epsilon} provides an analytical characterization of the $\epsilon$-subgradient step that may emerge in the considered setup; Lemma~\ref{lemma:bounded-error} establishes the constraints that~\eqref{eq:bounded_gradient_error} imposes on the tracking error $\|\by[t_k] - \bu[t_k]\|_2$; and finally, Theorem~\ref{thm:convergence} leverages Lemma~\ref{eq:epsilon} to establish asymptotic convergence of the reference signal $\bu[t_k]$ and the iterates $\bV[t_k]$ to the optimal solution of~$\mathrm{(P1)}$.\footnote{Proofs are omitted due to space constrains, and are available in~\cite{DhopleKKT}.}

\vspace{.1cm}

\begin{lemma}
\label{eq:epsilon}
If at time $t_k$, $\by_i[t_{k}] \neq \bu_i[t_{k}]$ for at least one dynamical system, then $\bh(\bV[t_k]) + \bg(\by[t_k], \bd)$ is an $\epsilon$-subgradient of the dual function at $\blambda[t_k]$.  In particular, under \emph{Assumption~\ref{ass:bounds_error}}, it holds that 
\begin{subequations}
\label{eq:eps_subgradient_all}
\begin{align}
& \hspace{-.3cm} \left(\bh(\bV[t_k]) + \bg(\by[t_k], \bd) \right)^\sfT (\blambda - \blambda[t_k])   \nonumber \\
& \hspace{2.5cm} \geq q(\blambda) - q(\blambda[t_k]) - \epsilon[t_k] \,\,\,\, \forall \,\, \blambda \, , \label{eq:eps_subgradient}
\end{align} 
where the error $\epsilon[t_k] \geq 0$ can be bounded as $\epsilon[t_k] \leq 2 \alpha_k \tilde{G} G^2$. \hfill $\Box$
\end{subequations}

\vspace{.1cm}

\emph{Proof.} For notational convenience, define 
\begin{subequations}
\begin{align}
\bs_u[t_k] &:= \bh(\bV[t_k]) + \bg(\bu[t_k], \bd), \nonumber \\
\bs_y[t_k] &:= \bh(\bV[t_k]) + \bg(\by[t_k], \bd).
\end{align}
Notice that $\bs_u[t_k]$ and $\bs_y[t_k]$ are the gradients of the dual function~\eqref{eq:dualFunction} evaluated at $\blambda[t_k]$ and $\tilde \blambda[t_k]$, respectively~\cite{Bertsekas_ConvexAnalysis}; i.e., it holds that 
\begin{align}
\bs_u^\sfT[t_k]  (\blambda - \blambda[t_k]) & \geq q(\blambda) - q(\blambda[t_k]),  \label{eq:graddefu}
\\
\bs_y^\sfT[t_k]  (\blambda - \tilde \blambda[t_k]) & \geq q(\blambda) - q(\tilde \blambda[t_k]), \forall \blambda \label{eq:graddef}
\end{align} 
Adding $\bs_y^\sfT[t_k] (\tilde \blambda[t_k] - \blambda[t_k])$ to both sides of~\eqref{eq:graddef}, one gets that the following holds $\forall \blambda$  
\begin{equation*}
\bs_y^\sfT[t_k]  (\blambda - \blambda[t_k]) \geq q(\blambda) - q(\tilde \blambda[t_k]) + \bs_y^\sfT[t_k] (\tilde \blambda[t_k] - \blambda[t_k]) \, .
\end{equation*}
Adding and subtracting $q(\blambda[t_k])$ to the right-hand-side of the inequality above, 
\begin{align}
&\bs_y^\sfT[t_k]  (\blambda - \blambda[t_k]) \geq \,\,q(\blambda) - q(\blambda[t_k])   \nonumber \\
& +q(\blambda[t_k]) - q(\tilde \blambda[t_k]) + \bs_y^\sfT[t_k](\tilde \blambda[t_k] - \blambda[t_k]). \label{eq:ineq}
\end{align}
With regard to~\eqref{eq:ineq}, define 
\begin{equation}
\epsilon[t_k]:= q(\tilde \blambda[t_k]) - q(\blambda[t_k]) + \bs_y^\sfT[t_k](\blambda[t_k] - \tilde \blambda[t_k]). \label{eq:edef}
\end{equation}
By using the definition of the gradient of the dual function at $\blambda[t_k]$ and applying the Cauchy-Schwartz inequality, one has that
\begin{align}
\epsilon[t_k] &\leq  \bs_u^\sfT[t_k](\tilde \blambda[t_k] - \blambda[t_k]) + \bs_y^\sfT[t_k](\blambda[t_k] - \tilde \blambda[t_k])\label{eq:edef0}  \\
&\leq  2 G \,\, \|\tilde \blambda[t_k] - \blambda[t_k]\|_2 \, \label{eq:edef1} \\
&\leq 2 \alpha_k \tilde G G^2 \label{eq:edef2}
\end{align}
\end{subequations}
where~\eqref{eq:bounded_gradient_2} was used to obtain~\eqref{eq:edef1} from~\eqref{eq:edef0}, and~\eqref{eq:edef2} follows from~\eqref{eq:step3}. \hfill $\Box$
\end{lemma}

Condition~\eqref{eq:bounded_gradient_error} implicitly bounds the tracking error $\|\by[t_k] - \bu[t_k]\|_2$, as specified in the following lemma. 

\vspace{.1cm}

\begin{lemma}
\label{lemma:bounded-error}
Under \emph{Assumption~\ref{ass:bounds_error}}, it follows that the tracking error $\|\by[t_k] - \bu[t_k]\|_2$, $k \in \mathbb{N}$, can be bounded as 
\begin{align}
\|\by[t_k] - \bu[t_k]\|_2 \leq \|-\bA^{-1} \bC^\sfT\|_2 \tilde{G} G \alpha_k. \label{eq:bounded_tracking_error}
\end{align}

\emph{Proof.} From the non-expansive property of the projection operator, the left-hand side of~\eqref{eq:bounded_tracking_error} can be bounded as:  
%
\begin{align}
\|\by[t_k] - \bu_i[t_k]\|_2 & \leq \|-\bA^{-1} \bC^\sfT(\tilde{\blambda}[t_k] - \blambda[t_k])\|_2 \nonumber \\
& \hspace{-0.4in} \leq \|-\bA^{-1} \bC^\sfT\|_2 \|\tilde{\blambda}[t_k] - \blambda[t_k]\|_2 \nonumber \\
& \hspace{-0.4in} \leq  \|-\bA^{-1} \bC^\sfT\|_2 \tilde{G} G \alpha_k  \label{eq:proof-track-s3} 
\end{align}
where~\eqref{eq:proof-track-s3} follows from~\eqref{eq:step3}.\hfill $\Box$ 
\end{lemma}

%

\vspace{.1cm} 
It can be noticed from~\eqref{eq:bounded_tracking_error} that the tracking error is allowed to be arbitrarily large, but 
the system output $\by[t_k]$ should eventually follow the reference signal $\bu[t_k]$. In fact, since the sequence $\{\alpha_k\}$ is majorized by $\{\eta_k\}$, and $\eta_k \downarrow 0$, it follows that $\|\by[t_k] - \bu[t_k]\|_2 \rightarrow 0$ as $k \rightarrow \infty$. 

While~\eqref{eq:bounded_tracking_error} bounds the error $\|\by[t_k] - \bu[t_k]\|_2$,   
asymptotic convergence of the reference signal $\bu[t_k]$ as well as of iterates $\bV[t_k]$ to the optimal solution of the steady-state optimization problem $\mathrm{(P1)}$ is established next. 
   
\vspace{.2cm}

\begin{theorem}
\label{thm:convergence}
Under \emph{Assumptions~\ref{ass:DynSystems}--\ref{ass:bounds_error}}, and using a stepsize sequence $\{\alpha_k\}_{k\geq 0}$ satisfying conditions $\mathrm{(s1)}$--$\mathrm{(s3)}$, the following holds for the closed-loop system~\eqref{eq:dualsubgradient_system}: 

\noindent \emph{i)} $\{\blambda[t_k]\}_{k \geq 0}$ is bounded; 

\noindent \emph{ii)} $\blambda[t_k] \rightarrow \blambda^{\mathrm{opt}}$ as $k \rightarrow \infty$;  

\noindent \emph{iii)} $\bV[t_k] \rightarrow \bV^{\mathrm{opt}}$  and $\{\bu_i[t_k] \rightarrow \bu_i^{\mathrm{opt}}\}$ as $k \rightarrow \infty$;  and, 

\noindent \emph{iv)} $\by_i(t) \rightarrow \bu_i^{\textrm{opt}}$ as $t \rightarrow \infty$, $\forall i \in \cN$. 

\noindent Statements \emph{i)}--\emph{iv)} hold for any initial conditions $\bV[0], \{\bu_i[0]\}, \{\by_i(0)\},\blambda[0]$, and $0 < t_{k} - t_{k-1} < \infty$, $k \in \mathbb{N}$. 

\vspace{.1cm}

\emph{Proof.} \emph{i)}--\emph{ii)} Boundedness and convergence of the dual iterates can be proved by leveraging the results of Theorem~3.4 in~\cite{Kiwiel04}. In particular, it suffices to show that the following technical requirement is satisfied in the present setup:  
\begin{subequations}
\begin{align}
\sum_{k = 0}^{+ \infty} \alpha_k \epsilon[t_k] <  + \infty \, . \label{proofT1-s1}
\end{align}
From Lemma~\ref{eq:epsilon}, it follows that the left-hand-side of~\eqref{proofT1-s1} can be bounded as 
\begin{align}
\sum_{k = 0}^{+ \infty} \alpha_k \epsilon[t_k] & \leq \sum_{k = 0}^{+ \infty}  2 \alpha_k^2 \tilde{G} G^2 \leq 2  \tilde{G} G^2 \sum_{k = 1}^{+ \infty}  \eta_k^2 \label{proofT1-s3}
\end{align}
\end{subequations}
where the second inequality in~\eqref{proofT1-s3} follows from the fact that $\alpha_k \leq \eta_k$ for all $k$. Since $\sum_{k = 0}^{+ \infty} \eta_k^2 < + \infty$, the series~\eqref{proofT1-s3} is finite, and thus~\eqref{proofT1-s1} holds.

\noindent \emph{iii)} From the strict convexity of the Lagrangian in the primal variables, it follows that optimal primal variables can be uniquely recovered as $\{ \bV^{\mathrm{opt}} \bu^{\textrm{opt}}\} = \arg \min_{\bV \in \cV, \bu \in \cU} L\left(\bV, \bu, \blambda^{\mathrm{opt}}\right)$. 

\noindent \emph{iv)} At convergence, the reference signal is constant, with value $\bu_i^{\textrm{opt}}$. Then, $\by_i(t) \rightarrow \bu_i^{\textrm{opt}}$ as $t \rightarrow \infty$ by~\eqref{eq:sys1-equilibrium}. \hfill $\Box$
\end{theorem}

\emph{Remark}. Problem~\eqref{eq:p1} could be solved either centrally or in a decentralized fashion, and the reference signals $\{\bu_i^{\textrm{opt}}\}_{i = 1}^N$ could be subsequently dispatched for the dynamical systems. It is evident that with these solutions the  optimization and local control tasks operate at two different time scales, with reference signals updated every time that problem~\eqref{eq:p1} is solved. Further, if relevant problem parameters change during the solution of~\eqref{eq:p1}, strategies operating under  time-scale separation would dispatch outdated setpoints. In contrast, steps~\eqref{eq:dualsubgradient_system} continuously pursue solutions of the formulated optimization problem by dynamically updating the setpoints, based on current system outputs and problem parameters.   \hfill $\Box$


\section{Dynamic Controller for Inverters}
\label{sec:ApplicationExample}

Consider a distribution system comprising $N+1$ nodes collected in the
set $\cN := \{0,1,\ldots,N\}$, with node $0$ denoting the secondary
of the transformer, and lines represented by the set of undirected 
edges $\cE := \{(m,n)\}$.
 For simplicity of exposition, assume that the system is balanced, and 
renewable-interfaced inverters are located at $\cN \backslash \{0\}$. However, the framework can be 
readily extended to account for: unbalanced multi-phase systems 
(by following the method in~\cite{Dallanese-TSG13}); load control 
(by considering four-quadrant inverters); and, nodes with no power generation 
(by adding relevant constraints~\cite{Dorfler14}).

Let $V_i\in \mathbb{C}$ and $I_i \in \mathbb{C}$ denote the phasors
for the line-to-ground voltage and the current injected at node $i \in
\cN$, respectively, and define $\bi := [I_0, \ldots, I_N]^\sfT \in
\mathbb{C}^{N+1}$ and $\bv := [V_0, \ldots, V_N]^\sfT \in
\mathbb{C}^{N+1}$. Using Ohm's and Kirchhoff's circuit laws, the linear
relationship $\bi = \bY \bv$ can be established, where the system
admittance matrix  $\bY \in \mathbb{C}^{N+1 \times N+1}$ is formed
based on the system topology and the $\pi$-equivalent circuits 
of the lines $(m,n) \in \cE$; see e.g.,~\cite{LavaeiLow, Dallanese-TSG13}.  

Similar to e.g.,~\cite{LavaeiLow, Dallanese-TSG13,OID,OID_TEC}, 
consider expressing powers and voltage magnitudes as linear functions of the
outer-product matrix $\bV := \bv \bv^\sfH$. Specifically, define 
$\mathbf{Y}_i := \be_i \be_i^\sfT \mathbf{Y}$ per node $i$, where $\{\mathbf{e}_{i}\}_{i \in \cN}$
denotes the canonical basis of $\mathbb{R}^{N+1}$. Based on
$\mathbf{Y}_i$, define also the Hermitian matrices $\bPhi_{i} :=
\frac{1}{2} (\bY_i + \bY_i^\sfH)$, $\bPsi_{i} := \frac{j}{2} (\bY_i -
\bY_i^\sfH) $, and $\bUpsilon_{i} := \be_i \be_i^\sfT$. Then, the net injected powers at node 
$i \in \cN$ can be  expressed as $\trace(\bPhi_i \bV)   =  \bar{P}_{i} - \bar{P}_{\ell,i}$ and 
$\trace(\bPsi_i \bV)   = \bar{Q}_{i} - \bar{Q}_{\ell,i}$, respectively, 
where $\bar{P}_{\ell,i}$ and $\bar{Q}_{\ell,i}$ denote the active and reactive
setpoints for the demand at node $i \in \cN \backslash \{0\}$, whereas $\bar{P}_{i}$ and $\bar{Q}_{i}$ are the 
active and reactive powers generated.  Further, $|V_i|^2$
is given by  $|V_i|^2 = \trace(\bUpsilon_i \bV)$. 

Upon denoting as $V_{\mathrm{min}}$ and $V_{\mathrm{max}}$ lower and upper limits 
for $\{|V_n|\}_{n = 0}^N$, matrix $\bV$ is confined to lie in the set  
\begin{align}
\cV^{1} :=  \{\bV \succeq \mathbf{0}: \rank(\bV) = 1, V_{\mathrm{min}}^2 \leq \trace(\bUpsilon_i \bV)  \leq V_{\mathrm{max}}^2 \,\forall \, i \} \, \nonumber  
\end{align} 
with the constraint $|V_0| = 1$ left implicit. 

For prevailing ambient conditions, let $P_i^{\textrm{av}}$ denote the
available active power for the inverter at node $i \in \cN \backslash \{0\}$. 
Then, the allowed operating regime for the inverter 
at node $i$ is assumed to be $\cY_i  =  \{ \bar{P}_{i}, \bar{Q}_{i} \hspace{-.1cm} :  0  \leq \bar{P}_{i}  \leq  P_{i}^{\textrm{av}}, \bar{Q}_{i}^2  \leq  S_{i}^2 - \bar{P}_{i}^2, 
| \bar{Q}_{i} | \leq \tan \theta \bar{P}_{i}\}$ (see e.g.,~\cite{OID,OID_TEC} for a more 
detailed explanation of possible inverter operating regions)
where $S_i$ is the apparent power rating, and $\theta$ models 
 power factor constraints. Set $\cY_i$ clearly adheres to \emph{Assumption~\ref{ass:compactSetY}}.

Powers $\{\bar{P}_{i},\bar{Q}_{i}\}_{i = 1}^N$ as well as the voltage-related matrix $\bV$ 
model the \emph{steady-state} operation of the distribution system.
For given load and ambient conditions, a prototypical OPF formulation for optimizing the
steady-state operation of the distribution system can be obtained by constraining 
variables $\bV$ and $(\bar{P}_{i}, \bar{Q}_{i})$ to the sets $\cV^{1}$ and $\cY_i$, respectively,  
and using the following mapping between the quantities explained above with the ones in~\eqref{eq:p1}:
\begin{subequations} 
\begin{align}
& \bu_i  = [\bar{P}_{i},\bar{Q}_{i}]^\sfT , \,\, \bd_i =  [\bar{P}_{\ell, i},\bar{Q}_{\ell, i}]^\sfT \, , \bC_i  = - \bI_{2 \times 2} ,\label{corr-1} \\
& \bD_i = \bI_{2 \times 2} , \,\,\, \bh_i(\bV)  = [\trace(\bPhi_i \bV), \trace(\bPsi_i \bV)]^\sfT  \label{corr-2}
\end{align} 
\end{subequations} 
and $H(\bV) =  \frac{a}{2} (\trace(\bPhi_0 \bV))^2 + b \trace(\bPhi_0 \bV)$, $a > 0,  b \geq 0 $. With this mapping,~\eqref{eq:p1-eq} represents the per-node balance equation for active and reactive powers, and $H(\bV)$ captures the cost of power 
drawn from (or supplied to) the substation. Unfortunately, the resultant 
optimization problem is nonconvex because of the constraint $\mathrm{rank}(\bV) = 1$. However, in the spirit of semidefinite relaxation~\cite{LavaeiLow}, this constraint can be
dropped; thus, $\cV^{1}$ is replaced by the convex set 
\begin{align}
\cV :=  \{\bV \succeq \mathbf{0}: V_{\mathrm{min}}^2 \leq \trace(\bUpsilon_i \bV)  \leq V_{\mathrm{max}}^2 \,\forall \, i \} \, .\label{eq:V_region}   
\end{align} 
Using $\cV$, problem $\mathrm{(P1)}$ turns out to be a relaxation of the AC OPF problem. If the optimal solution has
$\mathrm{rank}(\bV^{\mathrm{opt}}) = 1$, then the resultant power flows are globally optimal~\cite{LavaeiLow,Dallanese-TSG13}. 

\subsection{Dynamic Controller}
\label{sec:DynamicController}
 
Let $P_i(t)$ and $Q_i(t)$ denote the active and reactive powers of inverter $i$ averaged over one AC cycle, respectively, and consider setting state and output of the system~\eqref{eq:sys1} as follows: 
\begin{align}
\bx_i(t)  = [P_{i}(t),Q_{i}(t)]^\sfT \, , \,\, \, \, \by_i(t)  = \bx_i(t) \, .  \label{corr-inv} 
\end{align} 
Thus,~\eqref{eq:sys1-st} models the dynamics of real- and reactive-power controller at each inverter $i$, whereas~\eqref{eq:sys1-obs} boils down to a measurement of the inverter outputs. Dynamic models for the real and reactive power for inverters operating in a grid-connected mode are discussed in e.g.,~\cite[Ch.~8]{Iravanibook10}. Finally, variables $\bz_i(t)$ correspond to voltages on the lines connecting the inverters.  

The goal of the controller~\eqref{eq:dual_ascent_sys}--\eqref{eq:primal_y_sys}  is to steer the power output $\by_i(t) = [P_{i}(t),Q_{i}(t)]^\sfT$ of each inverter $i$ towards the OPF solution $\bu_i^{\mathrm{opt}} = [\bar{P}_{i}^{\mathrm{opt}},\bar{Q}_{i}^{\mathrm{opt}}]^\sfT$. During each interval $(t_{k-1}, t_{k}]$, the reference level $\bu_i[t_k] = [\bar{P}_{i}[t_k],\bar{Q}_{i}[t_k]]^\sfT$ is updated locally at each inverter $i$ via~\eqref{eq:primal_y_sys}, based on the most up-to-date multiplier $\blambda_i[t_k]$; while, the primal variable $\bV[t_k]$ is updated by a central authority (e.g.,  utility company), which aims to optimize the network performance. Thus, the resulting scheme is naturally \emph{decentralized}, and it entails a message passing that can be carried out via existing advance metering infrastructure protocols. Theorem 1 ensures that the output powers $P_{i}(t), Q_{i}(t)$ converge to $\bar{P}_{i}^{\mathrm{opt}},\bar{Q}_{i}^{\mathrm{opt}}$ for any slot duration $0 < t_{k} - t_{k-1} < \infty$, $k \in \mathbb{N}$; that is, in an effort to mitigate the strict time-scale separation between real-time inverter control and steady-state network optimization~\cite{Dorfler14}, steps~\eqref{eq:dual_ascent_sys}--\eqref{eq:primal_y_sys} are performed continuously (within the limits of affordable computational burden), without waiting for the inverters to converge to the intermediate reference levels, and without knowing the inverter controller dynamics.

\subsection{Representative numerical results}
\label{sec:NumericalResults}

Consider an illustrative low-voltage residential distribution system 
comprising a step-down transformer and $N= 5$ nodes featuring inverter-interfaced renewable sources. The node-node distance is set to $50$ m, and the line impedances are $Z_{mn} = 0.0135 + j 0.0045 \Omega$ for all $(m,n) \in \cE$. The optimization package \texttt{CVX} (\texttt{http://cvxr.com/cvx/}) is employed to perform the primal updates.  In the numerical test, the rank of matrix $\bV^\mathrm{opt}$ was $1$, meaning that the globally optimal solution of the OPF was identified~\cite{LavaeiLow}.  Voltage limits $V^\textrm{min}$ and $V^\textrm{max}$ are set to 0.95 pu and 1.05 pu, and the voltage magnitude at the substation is fixed to $|V_0| = 1$. The active and reactive loads are $1.10, 1.10, 1.10, 1.09, 1.10$ kW and $826, 828, 829, 821, 830$ VAr, respectively. As described in detail in~\cite{OID}, the operating regions $\{\cY_i\}$ of inverters providing ancillary services are formed based on the power ratings $\{S_i\}_{i = 1}^N$, which are taken to be $4.66, 4.83, 7.62, 7.62, 7.62$ kVA, as well as the available active powers $\{P_{i}^{\textrm{av}}\}_{i = 1}^N$, assumed to be $1.91, 1.95, 3.24, 3.24, 3.24$ kW.  In the cost function~\eqref{eq:p1-cost},  the cost of power drawn from the substation is $H(\bV) =  (\trace(\bPhi_0 \bV))^2 + 10 \times \trace(\bPhi_0 \bV)$, whereas $\bA_i = [1, 0; 0, .01]$ and $\bb_i = [10, 0.1]^\sfT$, for all $i = 1,\ldots, 5$. 

\begin{figure}[t]
\begin{center}
\subfigure[]{\includegraphics[width=8.8cm]{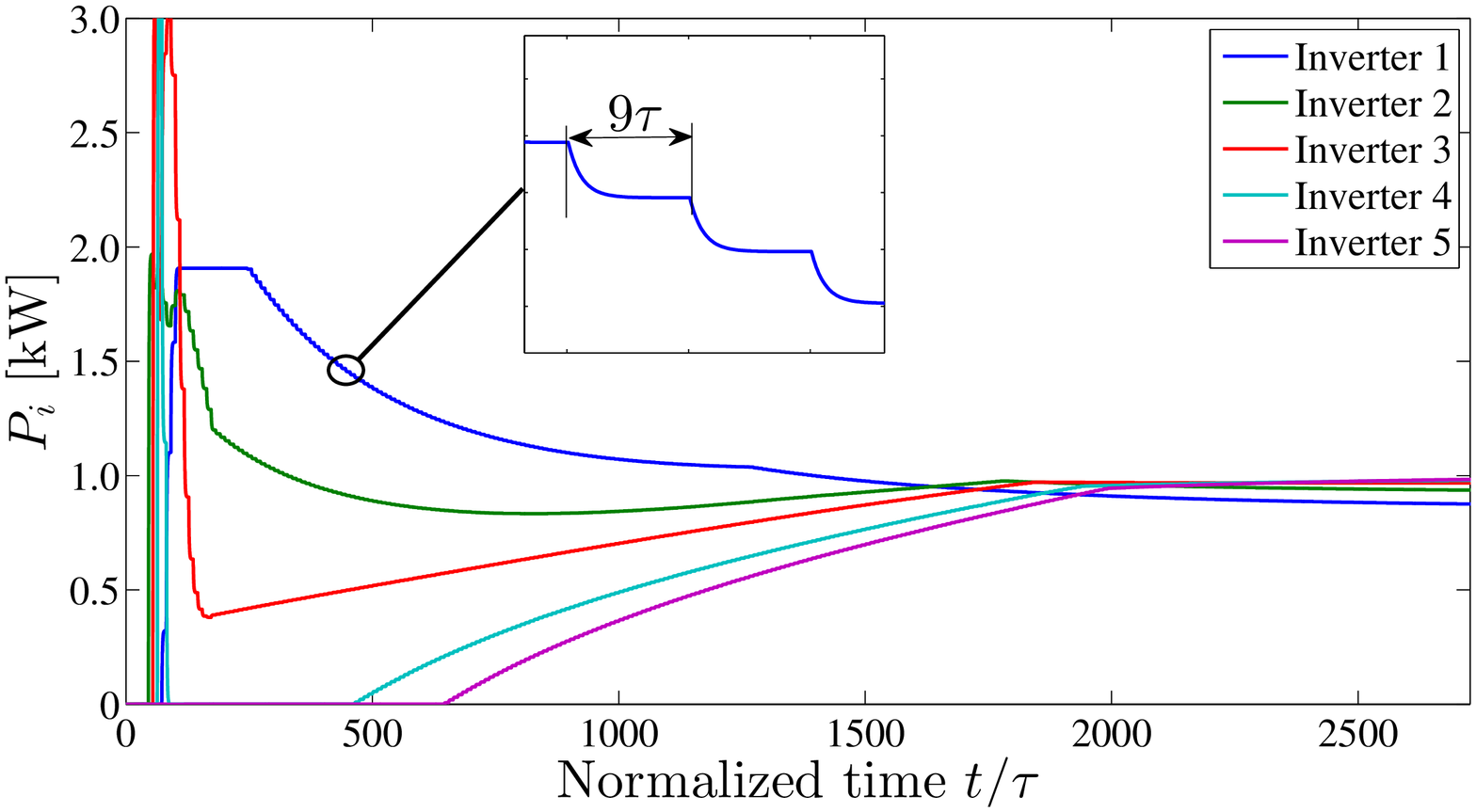}}
\subfigure[]{\includegraphics[width=8.8cm]{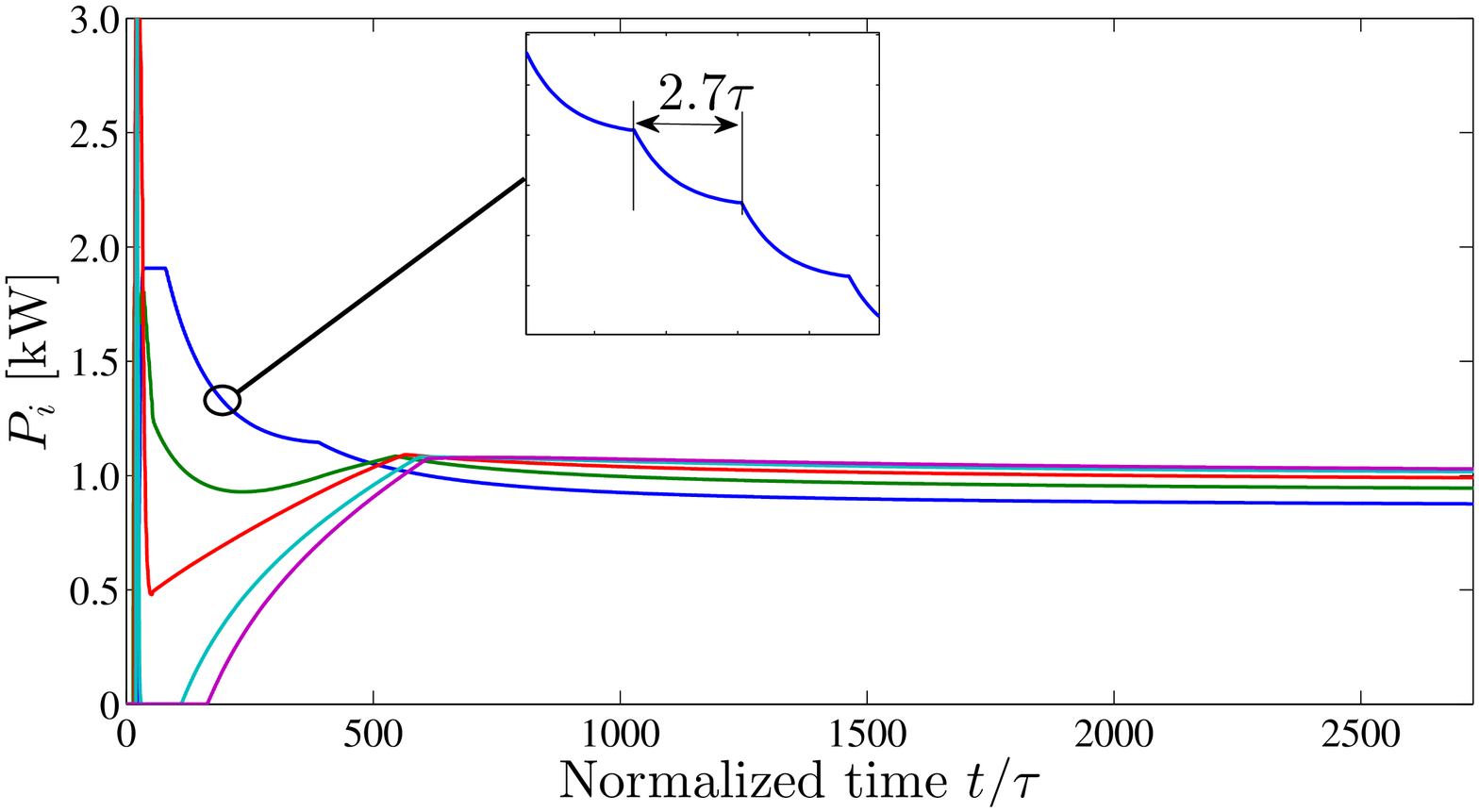}} 
\subfigure[]{\includegraphics[width=8.8cm]{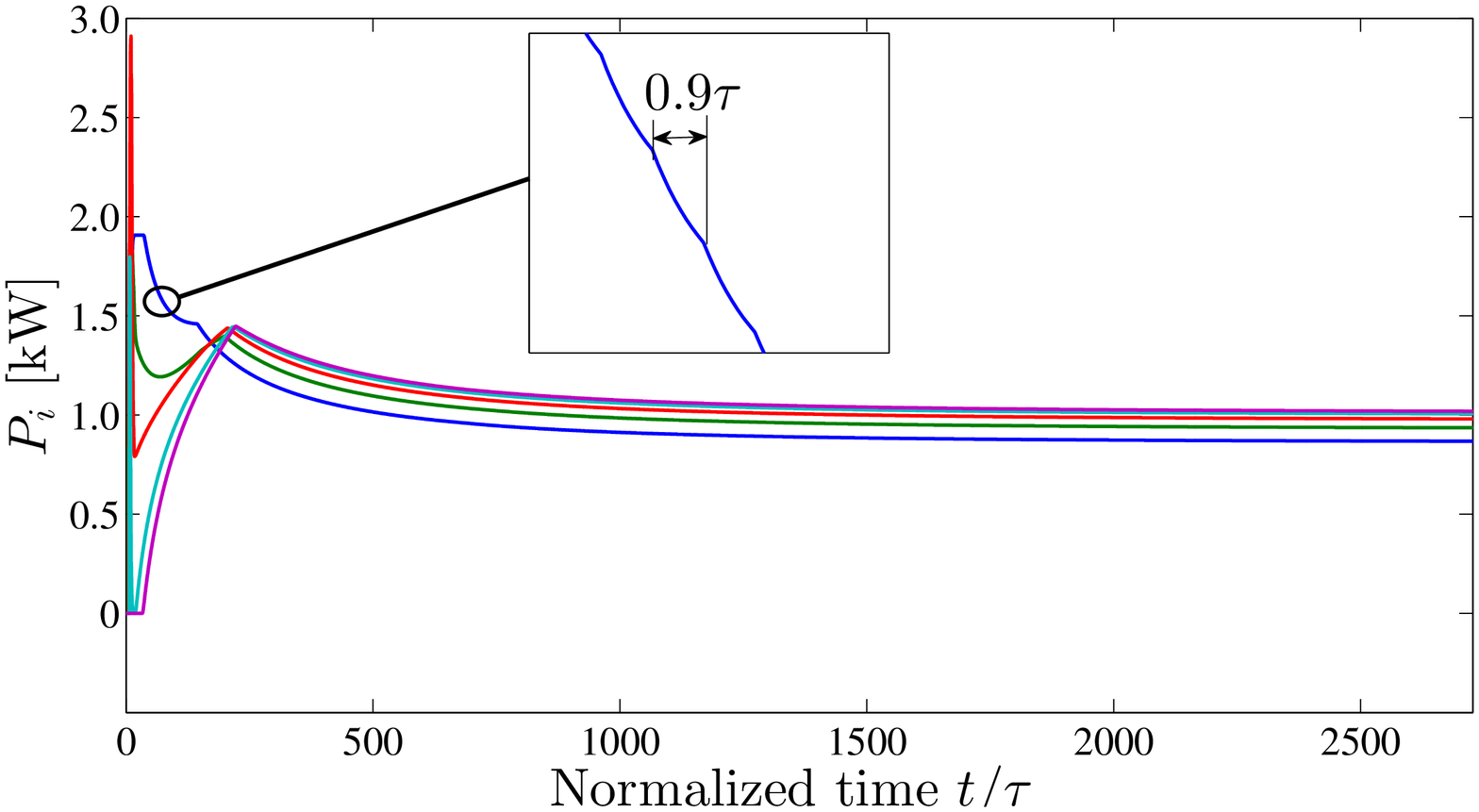}} 
\end{center}
\vspace{-.5cm}
\caption{Convergence of~\eqref{eq:dualsubgradient_system}, when the inverter-power dynamics are approximated as first-order systems with time constant $\tau$. Plots illustrate the evolution of the active powers when the system outputs are measured at intervals of length: (a) $9 \tau$; (b) $2.7 \tau$; and, $0.9 \tau$.  }
\label{Fig:convergence}
\vspace{-.5cm}
\end{figure}

A first-order system with time constant $\tau = 1.1$ is utilized to model the active and reactive power dynamics of each inverter~\cite[Ch.~8]{Iravanibook10}. For all $i = 1,\ldots, 5$, the initial states are set to $P_i(0) = Q_i(0) = 0$; for the voltage related matrix, the initial iterate is $\bV[0] = \bI$.  Finally, the stepsize $\alpha_k = 1/(10 \sqrt{k})$, for $k \geq 1$, is utilized. 

As a representative result, Fig.~\ref{Fig:convergence} illustrates the evolution of the active powers generated by the inverters (similar trajectories are obtained for the reactive powers, as well as for the active and reactive powers drawn from the point of common coupling).  Two setups are considered, depending on the duration of the interval $t_{k} - t_{k-1} $ between consecutive updates of variables $\bV[t_k]$ and $\bu[t_k]$: (a) $t_{k} - t_{k-1}  = 9 \tau$; (b) $2.7 \tau$, and, (b) $0.9 \tau$. Clearly, in setup (a), the output of the first-order systems converge to the reference inputs $\{\bu_i[t_k]\}$ within each interval $(t_{k-1}, t_{k}]$; see Fig.~\ref{Fig:convergence}(a). This yields a dual gradient step in~\eqref{eq:dual_ascent_sys} and, as expected, the overall scheme converges to the solution of the of the OPF $\mathrm{(P1)}$, which yields the following active powers: $0.86, 0.93, 0.97, 1.00, 1.01$ kW. 

In setup (b), updates of $\bV[t_k]$ and $\bu[t_k]$ are performed with higher frequency, and the output of the dynamical systems $\{\by_i[t_k]\}$ are different than the reference signals $\{\by_i[t_k]\}$ at each $t_k$, $k \geq 1$. Thus,~\eqref{eq:dual_ascent_sys} constitutes in this case an $\epsilon$-subgradient. Nevertheless,  outputs  $\{\by_i[t_k]\}$ converge to the solution of $\mathrm{(P1)}$, thus corroborating the claims of Theorem~\ref{thm:convergence}. In spite of the inexact update of the dual variables, the inverter outputs show a faster convergence to the solution of $\mathrm{(P1)}$ compared to Fig.~\ref{Fig:convergence}(a). In the third case (c), the time scales are further compressed, since the update of primal and dual variables is performed every $0.9 \tau$. By performing the primal and dual updates more frequently, this scheme converges markedly faster than (a) to the OPF solution.

\section{Concluding Remarks and Future Work}
\label{sec:Conclusions}

The paper developed a feedback controller for networked nonlinear dynamical systems, able to steer the system outputs to the solution of a convex constrained optimization problem. Global convergence was established for diminishing stepsize rules and strictly convex cost functions, even when the dynamical-system reference inputs are updated at a faster rate than the dynamical-system settling time.  The application of the proposed framework to the control of power-electronic inverters in AC distribution systems is discussed. 

\bibliographystyle{IEEEtran}
\bibliography{biblio.bib}
\end{document}